\documentclass[12pt,reqno,english,a4]{amsart}
\usepackage{amsmath,amsthm,amsfonts,amssymb,amscd}
\usepackage[latin1]{inputenc}
\usepackage[T1]{fontenc}
\usepackage[all]{xy}
\usepackage{psfrag}
\usepackage{epsfig}
\usepackage{a4wide}

\usepackage{hyperref} 
\theoremstyle{plain}
\newtheorem{theorem}{Theorem}
\newtheorem{proposition}{Proposition}
\newtheorem{lemma}{Lemma}

\newtheorem{remark}{Remark}
\newtheorem{corollary}{Corollary}
\newtheorem{definition}{Definition}

\numberwithin{equation}{section}

\newcommand{\om} {\omega}       
\newcommand{\op}{\operatorname}

\def \RR {{\mathbb R}}
\def \ZZ {{\mathbb Z}}

\def \TT {{\mathbb T}}
\def \SS {{\mathbb S}}
\def \cA {{\mathcal A}}
\def \cC {{\mathcal C}}
\def \cD {{\mathcal D}}
\def \cO {{\mathcal O}}
\def \cF {{\mathcal F}}
\def \hF {\widehat{\mathcal F}}
\def \cG {{\mathcal G}}

\def \cU {{\mathcal U}}

\def \cL {{\mathcal L}}

\def\-{{\setminus}}

\begin{document}

\title[On integrable  codimension one Anosov actions of $\RR^k$]
      {On integrable codimension one Anosov actions of $\RR^k$}

\author{Thierry Barbot}
\author{Carlos Maquera}
\thanks{The authors would like to thank FAPESP for the partial financial
support. Grants  2009/06328-2, 2009/13882-6 and 2008/02841-4 }



\date{\today}

\address{Thierry Barbot\\ Université d'Avignon et des pays de Vaucluse\\
LANLG, Faculté des Sciences\\
33 rue Louis Pasteur\\
84000 Avignon, France.}

\email{thierry.barbot@univ-avignon.fr}

\address{Carlos Maquera\\ Universidade de
S{\~a}o Paulo - S{\~a}o Carlos \\Instituto de ci{\^e}ncias
matem{\'a}ticas e de Computa\c{c}{\~a}o\\
Av. do Trabalhador S{\~a}o-Carlense 400 \\
13560-970 S{\~a}o Carlos, SP\\
Brazil}

 \email{cmaquera@icmc.usp.br}

 \maketitle
 \begin{abstract}
 In this paper, we consider codimension one Anosov actions of
 $\RR^k,\ k\geq 1,$ on closed connected orientable manifolds of dimension $n+k$
 with $n\geq 3$. We show that the fundamental group of the ambient manifold is solvable if and only if the weak foliation of codimension one is transversely affine. We also study the situation where one $1$-parameter subgroup of $\RR^k$ admits a cross-section, and compare this to the case where the whole action is transverse to a fibration over a manifold of dimension $n$.
 As a byproduct, generalizing a Theorem by Ghys in the case $k=1$, we show that, under some assumptions about the
 smoothness of the sub-bundle $E^{ss}\oplus E^{uu}$, and in the case where the action preserves the volume,
 it is
 topologically equivalent to a suspension of a linear Anosov action of
 $\mathbb{Z}^k$ on $\TT^{n}$.
 \end{abstract}

 \medskip

 \thispagestyle{empty}


 \section{Introduction}
 We consider the problem of topological classification of higher-rank
 abelian Anosov actions. More precisely, these are actions of $G=\ZZ^k$ or $\RR^k,
 \ k\geq2,$ by commuting diffeomorphisms of a compact manifold with at least one Anosov element,
 that is, some element $a\in \RR^k$ acts normally hyperbolically with respect to the orbit
 foliation.

   This concept was originally introduced by Pugh--Shub \cite{push} in the early
 seventies and, when $k=1$, coincides with the classical notion of Anosov
 diffeomorphism (in the discrete case) and Anosov flow.
 Recently this actions received a strong impetus under the contribution
 of A. Katok and R.J. Spatzier. The rigidity aspects of these actions
 received in the last decade a lot of attention, in the framework of Zimmer program.

 Let $\phi$ be an Anosov action of $\RR^k$ and $a\in \RR^k$ such that
 $g:=\phi(a,\cdot)$ is normally hyperbolically with respect to the orbit
 foliation. As for Anosov flows, there exists a continuous $Dg$-invariant  splitting of the tangent
 bundle
  $$
  TM=E^{ss}\oplus T\cO \oplus E^{uu}
  $$ such that
  $$
  \begin{array}{cc}
   \|Dg^n|_{E^{ss}}\|\leq Ce^{-\lambda n} & \forall n>0 \\
   \|Dg^n|_{E^{uu}}\|\leq Ce^{\lambda n} & \forall n<0
  \end{array},
  $$
  Where $T\cO$ denotes  the $k$-dimensional subbundle of $TM$ that is tangent to the orbits of 
  $\phi.$

   Hirsch, Pugh and Shub developed the basic theory of normally hyperbolic transformations
 in \cite{hirpush}. As consequence of this we obtain that the splitting is H\"older continuous
 and the subbundles $E^{ss}$, $E^{uu}$, $T\cO \oplus E^{ss}$, $T\cO \oplus E^{uu}$
 are integrable.
 The corresponding foliations, $\cF^{ss}$, $\cF^{uu}$, $\cF^{s}$, $\cF^{u}$, are called
 \textit{the strong stable foliation, the strong unstable foliation, the weak stable foliation},
 and \textit{the weak unstable foliation}, respectively.

 If the foliation $\cF^{uu}$ is one dimensional, we say that $\phi$ is of \textit{codimension
 one}. We note that when an Anosov action $\phi$ of $\ZZ^k$ is of codimension
 one, then, by the well-know result of Franks-Newhouse (\cite{franks}, \cite{newhouse}), the Anosov element is known
 to be topologically conjugated to a hyperbolic toral automorphism. This
 gives topological conjugacy of the whole action to an Anosov linear
 action of $\ZZ^k$ on the torus.

 For Anosov flows, in 1970's Verjovsky conjectured that: \textit{every codimension one Anosov flow
 on a manifold of dimension $\geq 4$ is topologically equivalent to a suspension (of a hyperbolic toral
 automorphism)}. The hypothesis on the dimension is necessary, since geodesic flows of negatively curved surfaces provide counter-examples in dimension $3$.
 
 It is a folkloric fact that to prove Verjovsky's conjecture amounts to find a cross-section to the Anosov flow, i.e. a codimension one closed submanifold transverse to the flow intersecting every positive orbit of the flow: indeed, if this occurs, the flow is obviously topologically
 equivalent to the suspension of the first return map, which itself is clearly a codimension one Anosov diffeomorphism. Franks-Newhouse Theorem stated above provide then the conclusion.

 Results on the existence of global cross section for codimension one Anosov flows on a manifold of dimension $\geq 4$ have been obtained by many authors in different contexts (see, \cite{plante1}, \cite{plante2}, \cite{ghys} and \cite{simic}).
 Many of them are based on the Schwartzman criterion (see \S~\ref{sub:schwartzman}). For example, J.F. Plante used it to show that Anosov flows for which the weak stable foliation is transversely affine admits a cross-section (\cite{plante2}). It allows to prove that the same conclusion holds if the fundamental group of the ambient manifold is solvable (see \cite{matsumotosolv} for a complement to Plante's argument in higher dimensions). Let us also mention Ghys Theorem (\cite{ghys}): if $E^{ss}\oplus E^{uu}$ is of class
 $C^1$ and $\dim M > 3$ then the flow admits a global cross section.

 In this paper we bring the study to the case of codimension one Anosov
 actions of $\RR^k, k>1$ (see Definition~\ref{defi:defcodone}). For example, we will prove, generalizing Ghys Theorem stated above:

 \begin{theorem}
 \label{thm:main}
 Let $\phi$ be a $C^2$ codimension one Anosov action of $\RR^k$ on a closed manifold
 $M$ of dimension $k+n$ with $n\geq 3$. If $E^{ss}\oplus E^{uu}$ is $C^1$ and if $\phi$ is volume-preserving, then $\phi$ is topologically equivalent to a suspension of
 a linear Anosov action of $\mathbb{Z}^k$ on $\TT^{n}$.
 \end{theorem}
 
 More generally, it seems reasonable to conjecture, as already mentioned in \cite{Bar-Maq}\footnote{Let's make precise that Verjovsky didn't formulate himself this conjecture in this more general form.}:
 
  \vspace{.5cm}
 \noindent
 \textbf{Verjovsky's conjecture for $\RR^k$-actions.}
 \textit{Every irreducible codimension one Anosov action of $\RR^k$ on a manifold of dimension at least $k+3$ is topologically conjugate to the suspension of a linear Anosov action of $\ZZ^k$ on the torus. }

 \vspace{.5cm}

Actually, it has also been conjectured that such a topological conjugacy, in the case $k \geq 2$, should be smooth (see for example \cite[Conjecture 1.1]{kalininspatz}), hence we can replace in this conjecture the term <<topologically conjugate>> by <<smoothly conjugate>>\footnote{Observe that  <<irreducible>> has not exactly the same meaning in \cite{kalininspatz} and in  \cite{Bar-Maq}; however, irreducible codimension actions in the meaning of Verjovsky conjecture above are irreducible in the meaning of \cite{kalininspatz}. }(under the hypothesis $k\geq 2$).

However, in the case $k\geq 2$, we have to face new difficulties, even in the <<simple cases>> mentioned above (for example, when the fundamental group of the manifold $M$ is solvable), which are trivially solved in the case $k=1$ but not anymore in the general case. 

The arguments mentioned above in the context $k=1$, like Schwartzman criterion, leads naturally, in the case $k\geq 2$, to the existence of a cross-section to one flow corresponding to a one parameter subgroup of $\RR^k$: in that situation, it means that the action is topologically equivalent to a \textit{$1$-suspended action,} ie. the suspension of a partially hyperbolic diffeomorphism $f$ whose central direction is tangent to a locally free action of $\RR^{k-1}$ commuting with $f$ (cf. Definition~\ref{def:susf}; hence, the generalized Verjovsky conjecture requires a classification of such partially hyperbolic diffeomorphisms). Whereas in the case $k=1$, the existence of a cross-section to the flow is equivalent to the existence of a fibration $\pi: M \rightarrow \SS^1$ whose fibers are transverse to the flow, it is far from obvious that $1$-suspended actions should be transverse to a fibration $\pi: M \rightarrow \TT^k$! Worse, even assuming that the action is transverse to such a fibration over $\TT^k$, and thus topologically equivalent to the suspension of an action of $\ZZ^k$ on the fibers (the monodromy representation, see Corollary~\ref{cor:susmonodromy}), we don't know how to prove that this action of $\ZZ^k$ is Anosov!

Now that we have listed what we \textit{should} prove, let us state what we prove in this paper. First of all, we will show that Plante's arguments can be generalized, providing the following Theorem:

\begin{theorem}
\label{teo:transaffine}
Let $\phi$ be a codimension one Anosov action on a closed manifold $M$. Then the following properties are equivalent:
\begin{enumerate}
\item the fundamental group of $M$ is solvable,
\item the weak stable foliation $\cF^s$ admits a transverse affine structure,
\item $\phi$ is splitting.
\end{enumerate}
\end{theorem}

In this statement, \textit{splitting} means that at the universal covering level, every strong stable leaf intersects every unstable leaf. This notion has proved to be an important one in the context of Anosov flows or diffeomorphisms: besides Theorem~\ref{teo:transaffine} itself, it is a crucial ingredient of Franks's study of codimension one Anosov diffeomorphisms (\cite{franks}). 

We will actually prove a stronger version (see Theorem~\ref{teo:transaffinecomplet}): we will show that the transverse affine structure, if it exists, is automatically $C^1$, and is \textit{complete} (in particular, the developping map is an diffeomorphism onto $\RR$). 

From now on, codimension one actions satisfying the equivalent assertions of Theorem~\ref{teo:transaffine} are called \textbf{splitting}.

Using Schwartzman criterion we will prove that \textit{splitting Anosov actions are topologically equivalent to $1$-suspended actions} (cf. Theorem~\ref{teo:splitsusp}).

Then, in \S~\ref{sec:fibration} we will turn our attention to the case where $\phi$ is transverse to a fibration $\pi: M \rightarrow B$. We show that up to a finite covering, $B$ is diffeomorphic to the torus $\TT^k$. We will collect several facts showing that the action have many properties similar to suspensions of toral automorphisms. In particular:

\begin{itemize}
\item the action is splitting,
\item the fibers of the fibration are homotopy equivalent to $\TT^n$ (hence, by a C.T.C. Wall Theorem homeomorphic to $\TT^n$, except maybe if $n=4$),
\item if one monodromy element is Anosov, then the action is topologically equivalent to the suspension of an action of $\ZZ^k$ by toral automorphisms.
\end{itemize}

As an application of all this study, we will prove Theorem~\ref{thm:main} stated above.
In the final section \S~\ref{sec:conclusion} we give some prospective about the way this work could be pursued further.

\section{Preliminaries}
\label{sec.preli}
 We start by establishing some definitions and notations. We consider an
 action $\phi :\RR^k \times M \to M$ of $\RR^k$ on a manifold $M$ of dimension $n+k.$
 We denote by
 $\mathcal{O}_p:=\{\phi(a, p),\om \in \RR^k \}$ the orbit of
 $p \in M$ and by $\Gamma_p := \{ a \in\RR^k : \phi(a, p) = p \}$ the isotropy group of $p.$ The action $\phi$ can be thought as
 a morphism $\RR^k \rightarrow \mbox{Diff}(M)$, the kernel of which is called
  the \textit{kernel of $\phi$:} the kernel is the intersection of all the isotropy groups.
 The action can also be thought as a Lie algebra morphism ${\zeta}: {\mathcal R}^k \rightarrow {\mathcal X}(M)$, where
 ${\mathcal R}^k$ denotes the Lie algebra of $\RR^k$ and ${\mathcal X}(M)$ the Lie algebra of vector fields on $M$.
 
We say that  $\phi $ and $\psi$, two actions of $\RR^k$ on a manifold $M$,  are {\it topologically equivalent}, if there is a homeomorphism $h:M\to M$ mappings orbits of $\phi$ onto orbits of $\psi$. 

  The action
 is said to be \textit{locally free} if the isotropy group of every point is
 discrete. Equivalently, it means that for every $\alpha$ in ${\mathcal R}^k\setminus \{0\},$ $\zeta(\alpha)$ is a
 non-vanishing vector field. In this case the orbits are diffeomorphic to
 $\mathbb{R}^{\ell} \times \mathbb{T}^{k-\ell}$, where $0\leq \ell \leq k$.
 For every $a$ in $\RR^k$ (respectively every $\alpha$ in ${\mathcal R}^k$), we denote by $\phi^a$ the diffeomorphism $x \mapsto \phi(a, x)$
 (respectively $X_\alpha$ the vector field $\zeta(\alpha)$).
 The flow $(t,x) \mapsto \phi(ta, x)$ will be called \textit{the flow generated by $a$,} and will be denoted by $\phi^{ta}$ (the context will help to avoid any confusion with the diffeomorphism $\phi^{ta}$ for a given real number $t$).

Let $\mathcal{F}$ be a continuous foliation on a manifold $M^{n+k}$. We
denote the leaf that contains $p \in M$ by $\cF(p)$.
We denote the tangent bundle of $M$ by $TM$. If $\cF$ is a $C^1$
foliation, then we denote the tangent bundle
 of $\cF$ by $T\cF$. The orbits of $\phi$ are the leaves
 of a central foliation $\cO$ that we call \textit{orbit foliation}. In particular, we denote by $T\cO$ the $k$-dimensional subbundle of $TM$
 that is tangent to the orbits of $\phi$.

 We fix a Riemannian metric $\|\cdot \|$, and denote by
 $d$ the associated distance map in $M$.

 \subsection{Anosov $\RR^k$-actions}
 \label{sub.def}

 \begin{definition}
 \label{defi:defcodone}
 {\rm
 Let $M$ be a $C^\infty$ manifold of dimension $n+k$ and $\phi$ a locally free $C^{2}$ action of
 $\RR^k$ on $M$.
 \begin{enumerate}
   \item We say that $a\in \RR^k$  is an \textit{Anosov element} for $\phi$ if $g=\phi^a$
   acts normally hyperbolically with respect to the orbit foliation. That is, there exist real
  numbers $\lambda > 0,\ C > 0$ and a continuous $Dg$-invariant splitting of the tangent bundle
  $$
  TM=E_a^{ss}\oplus T\cO\oplus E_a^{uu}
  $$
  such that
  $$
  \begin{array}{cc}
   \|Dg^n|_{E_a^{ss}}\|\leq Ce^{-\lambda n} & \forall n>0 \\
   \|Dg^n|_{E_a^{uu}}\|\leq Ce^{\lambda n} & \forall n<0
  \end{array}
  $$
   \item Call $\phi$ an \textit{Anosov action} if some $a\in \mathbb{R}^k$ is
   an Anosov element for $\phi$.
  \end{enumerate}

 The action $\phi$ is a {\it codimension-one} Anosov action
 if $E^{uu}_a$ is one-dimensional for some $a$ in $\mathbb{R}^k$.
 }
 \end{definition}

 Hirsch, Pugh and Shub developed the basic theory of normally hyperbolic transformations
 in \cite{hirpush}. As a consequence of this we obtain that the splitting is H\"older continuous
 and the subbundles $E_a^{ss}$, $E_a^{uu}$, $T\cO \oplus E_a^{ss}$, $T\cO \oplus E_a^{uu}$
 are integrable.
 The corresponding foliations, $\cF^{ss}[a]$, $\cF^{uu}[a]$, $\cF^{s}[a]$, $\cF^{u}[a]$, are called
 \textit{the strong stable foliation, the strong unstable foliation, the weak stable foliation},
 and \textit{the weak unstable foliation}, respectively.

 From now we will assume that every Anosov action comes with a specified Anosov element $a_0$.
 When $\phi$ is codimension one,
 we will always assume that $a_0$ has one dimensional strong unstable foliation.
 We denote by $\cF^u$, $\cF^s$, $\cF^{ss}$ and $\cF^{uu}$ the foliations associated to $a_0$.

 For all $\delta >0$, $\cF_{\delta}^{i}(x)$
 denotes the open ball in $\cF^i(x)$ centered at $x$ with radius
 $\delta$ with respect to the restriction of $\|\cdot \|$ to $\cF^i(x)$, where $i=ss,uu,s,u.$

 \begin{theorem}[of product neighborhoods]
 \label{thm:local product}
 Let $\phi:\RR^k \times M \to M$ be an Anosov action. There exists a
 $\delta_0 >0$ such that for all $\delta \in (0,\delta_0)$ and for
 all $x\in M,$ the applications
 $$
 [\cdot,\cdot]^u:\cF^s_\delta(x) \times \cF^{uu}_\delta(x)\to M; \ \
 [y,z]^u=\cF^s_{2\delta}(z)\cap \cF^{uu}_{2\delta}(y)
 $$
 $$
 [\cdot,\cdot]^s:\cF^{ss}_\delta(x)\times \cF^{u}_\delta(x)\to M; \ \
 [y,z]^s=\cF^{ss}_{2\delta}(z)\cap \cF^{u}_{2\delta}(y)
 $$

 are homeomorphisms onto their images.
 \end{theorem}

 \begin{remark}
 \label{rk.zero}
 {\rm
 Every foliation $\cF^{ss}$, $\cF^{uu}$, $\cF^{s}$ or  $\cF^{u}$ is preserved by every diffeomorphism
 commuting with $\phi^{a_0}$. In particular, it is
 $\RR^k$-invariant.  Another straightforward observation is that, since every compact domain in a leaf
 of $\cF^{ss}$ (respectively
 of $\cF^{uu}$) shrinks to a point under positive (respectively negative) iteration by $\phi^{a_0}$, every leaf
 of $\cF^{ss}$ or $\cF^{uu}$ is a \textit{plane}, i.e. diffeomorphic to $\RR^\ell$ for some $\ell$.

 Let $F$ be a weak leaf, let us say a weak stable leaf. For every strong stable leaf $L$ in $F$, let
 $\Gamma_L$ be the subgroup of $\RR^k$ comprising elements $a$ such that $\phi^a(L)=L$,
 and let $\omega_L$ be the saturation of $L$ under $\phi$. Thanks to Theorem~\ref{thm:local product} we have:
 \begin{itemize}
 \item $\omega_L$ is open in $F$,
 \item $\Gamma_L$ is discrete.
 \end{itemize}
 Since $F$ is connected, the first item implies $F=\omega_L$: the $\phi$-saturation of a strong leaf  is an entire weak leaf.
 Therefore, $\Gamma_L$ does not depend on $L$, only on $F$.
 The second item implies that the quotient $P=\Gamma_L\backslash\RR^k$ is a manifold, more precisely,
 a flat cylinder, diffeomorphic to $\RR^p\times\TT^q$ for some $p$, $q$.
 For every $x$ in $F$, define $\op{p}_F(x)$ as the equivalence class $a+\Gamma_L$ such that
 $x$ belongs to $\phi^a(L)$.  The map $\op{p}_F: F \to P$ is a locally trivial fibration
 and the restriction of $\op{p}_F$ to any $\phi$-orbit in $F$ is a covering map.
 Since the fibers are contractible (they are leaves of $\cF^{ss}$, hence planes), the fundamental group
 of $F$ is the fundamental group of $P$, i.e. $\Gamma_L$ for any strong stable leaf $L$ inside $F$.

 Observe that if $\cF^{ss}$ is oriented, then the fibration $\op{p}_F$ is trivial: in particular, $F$ is diffeomorphic to
 $P \times \RR^{p}$, where $p$ is the dimension of $\cF^{ss}$.

 Of course, analogous statements hold for the strong and weak unstable leaves.
 }
 \end{remark}

 \begin{remark}
 \label{rk.first}
 {\rm
 Let $\mathcal{A}=\mathcal{A(\phi)}$ be the set of Anosov elements of $\phi$.
  \begin{enumerate}
   \item $\mathcal{A}$ is always an open subset of $\mathbb{R}^k$.

   \item Every connected component of $\mathcal{A}$ is an open convex cone in
   $\mathbb{R}^k$.
   \end{enumerate}
  }
 \end{remark}

 For more details about this remark see \cite{Bar-Maq}.
 We call a connected component of $\mathcal{A}$ as a \textit{chamber}, by analogy with the case of
   Cartan actions. We denote by $   \cA_0$ the chamber containing the preferred Anosov element $a_0$.

   \subsection{Irreducible codimension one Anosov actions}
   \label{sub:reduc}
 A codimension one Anosov action $\phi$ of $\RR^k$ on $M$ is said to be
 \textit{irreducible} if for any $a\in \RR^k-\{0\}$ and $x\in M$
 with $\phi^a(x)=x$ we have that ${\rm Hol}_{\gamma}$, the holonomy along
 of $\gamma=\{\phi^{sa}(x);\ s\in [0,1]\}$ of $\cF^{s}(x)$, is
 a topological contraction or a topological expansion.

 \begin{theorem}[Theorem 7, \cite{Bar-Maq}]
 \label{thm:ReducingAction}
 Let $\phi:\RR^k\times M\to M$ be a codimension one Anosov action.
 Then, there exists a free abelian subgroup $H_0 \approx \RR^\ell$ of  $\RR^k$, a lattice
 $\Gamma_0 \subset H_0$, a smooth
 $(n+k-\ell)$-manifold $\bar{M}$, and $p:M\to \bar{M}$ a smooth
 $\TT^{\ell}$-principal bundle such that:
 \begin{enumerate}
   \item $\Gamma_0$ is the kernel of the action $\phi;$
   \item every orbit of ${\phi}_0=\phi|_{H_0 \times M}$ is a fiber of $p:M\to \bar{M}$.
   In particular, $\bar{M}$ is the orbit space of $\phi_0$;
   \item $\phi$ induces an irreducible codimension one Anosov action
   $\bar{\phi}: \bar{H} \times \bar{M} \to \bar{M}$ where $\bar{H}=\RR^k/H_0$.
    \end{enumerate}
 \end{theorem}

 In the terminology of \cite{brintop, bringroup}, Theorem~\ref{thm:ReducingAction} states
 that codimension one Anosov actions are $\TT^\ell$-extensions of irreducible codimension one Anosov actions.

 \begin{remark}
 \label{rk:solvablered}
 In particular, there is an exact sequence:
 $$\pi_2(\bar{M}) \rightarrow  \pi_1(\TT^\ell) \rightarrow \pi_1(M) \rightarrow \pi_1(\bar{M}) \rightarrow 0$$
 According to item $(3)$ above and Remark \ref{rk:redcase} the second homotopy group $\pi_2(\bar{M})$ vanishes.
 Hence $\pi_1(M)$ is an extension of $\pi_1(\bar{M})$ by $\ZZ^k$. In particular, $\pi_1(M)$ is solvable if and only if
 $\pi_1(\bar{M})$ is solvable.
 \end{remark}

  \subsection{The orbit space of a codimension one Anosov action}

 Let $\op{p}_{\op{univ}}:\widetilde{M}\to M$ be the universal covering map of $M$ and $\widetilde{\phi}$ be the
 lift of $\phi$ on $\widetilde{M}$. The foliations $\cF^{ss},\ \cF^{uu},\ \cF^{s}$ and
 $\cF^{u}$ lift to foliations $\widetilde{\cF}^{ss},\ \widetilde{\cF}^{uu},\ \widetilde{\cF}^{s}$ and
 $\widetilde{\cF}^{u}$ in $\widetilde{M}$. We denote by $Q^{\phi}$ be the
 orbit space of $\tilde{\phi}$ and $\pi^{\phi}:\widetilde{M}\to Q^{\phi}$ be the
 canonical projection. If $\phi$ is an irreducible codimension one Anosov action of $\RR^k$ on $M$, then
 the following properties can be found in Subsection 4.3 of \cite{Bar-Maq}:

 \begin{enumerate}
   \item $Q^{\phi}$, the orbit space of $\widetilde{\phi}$, is homeomorphic to $\RR^{n}.$

   \item The foliations $\widetilde{\cF}^{uu}$, $\widetilde{\cF}^{ss}$, $\widetilde{\cF}^{u}$,
  $\widetilde{\cF}^{s}$ and the foliation defined by $\widetilde{\phi}$ are by closed planes.
  The intersection between a leaf of $\widetilde{\cF}^{u}$ and a leaf of
  $\widetilde{\cF}^{s}$ is at most an orbit of $\widetilde{\phi}$.
  Every orbit of $\widetilde{\phi}$ meets a leaf of $\widetilde{\cF}^{uu}$ or $\widetilde{\cF}^{ss}$
  at most once.

   \item The universal covering of $M$ is diffeomorphic to $\RR^{n+k}.$

   \item All the non compact orbits are planes.
 \end{enumerate}

 \begin{remark}
 \label{rk:redcase}
 It is easily proved thanks to Theorem \ref{thm:ReducingAction} that all these statements, except $(4)$,
 are true even for reducible codimension one actions.
 \end{remark}

The fundamental group $\pi_1(M)$ acts on $\widetilde{M}$ by deck transformations. This action preserves the orbits of $\tilde{\phi}$. Hence, it induces an action of $\pi_1(M)$ by diffeomorphisms. Since every orbit of an irreducible Anosov action of codimension one is either an incompressible compact torus or a plane:

 \begin{lemma}
 Assume that $\phi$ is irreducible. Then, the stabilizer in $\pi_1(M)$ of every element $x$ of $Q^\phi$ is either trivial, or isomorphic to $\ZZ^k$.
 \end{lemma}

 \subsection{The leaf spaces associated to a codimension one Anosov action}

\begin{definition}
{\rm
The \textit{stable leaf space} is the quotient of $\widetilde{M}$ be the equivalence relation identifying two points if they lie in the same leaf of $\widetilde{\cF}^{s}$. We denote it by
$\cL^s$, and by $\pi^s: \widetilde{M} \rightarrow \cL^s$ the quotient map.
}
\end{definition}

The stable leaf space is an $1$-manifold: for example, the restrictions of $\pi^s$ to strong unstable leaves of $\widetilde{\mathcal F}^{uu}(\tilde{p})$ define an atlas of open subsets homeomorphic to
$\RR$. Observe that this atlas is $C^1$. Moreover, since every weak stable leaf is a closed plane, it disconnects $\widetilde{M}$. It follows that $\cL^s$ is simply connected.

However, $\cL^s$ may not satisfy the Hausdorff separation property, and hence is not necessarily diffeomorphic or homeomorphic to $\RR$.

For more details, see for example \cite{haereeb}, \cite{Bartree} or \cite{matsumotosolv}.

\begin{definition}
{\rm
The action $\phi$ is \textit{splitting} if every leaf of $\widetilde{\mathcal F}^{uu}$ intersects every leaf of $\widetilde{\mathcal F}^s$.
}
\end{definition}

This notion appeared in the context of Anosov diffeomorphisms: one crucial step in the classification Theorem of Anosov diffeomorphisms of codimension one is to prove that they are splitting (cf. \cite{franks}). Observe that in the splitting case, the restriction of $\pi^s$ to any strong unstable leaf realizes a diffeomorphism between $\cL^s$ and $\RR$.

Finally, since $\widetilde{\mathcal F}$ is tangent to the orbits, the projections by $\pi^\phi$ of two different leaves of $\widetilde{\mathcal F}^s$ (respectively $\widetilde{\mathcal F}^u$) are disjoints. These projections define two foliations, denoted respectively by $\cG^s$, $\cG^u$, of $Q^\phi \approx \RR^n$. Observe that the leaves of $\cG^u$ are closed lines, and the leaves of $\cG^s$ are closed hyperplanes diffeomorphic to $\RR^{n-1}$. 

Finally, $\pi^\phi$ realizes an diffeomorphism between the leaf space of $\cG^s$ and the leaf space $\cL^s$ of $\widetilde{\mathcal F}^s$.

Since the action of $\pi_1(M)$ by deck transformations preserves stables leaves, it induces an action of $\pi_1(M)$ on $\cL^s$.
The following proposition is an immediate corollary of the description of the leaves of $\cF^s$ of irreducible Anosov actions of codimension one: non-planar stable leaves are stable leaves of compact orbits.

\begin{proposition}
Assume that $\phi$ is irreducible. Then, the stabilizer in $\pi_1(M)$ of every element of $\cL^s$ is either trivial, or isomorphic to $\ZZ^k$.
\end{proposition}

Moreover, since a given stable or unstable leaf contains at most one compact $\phi$-orbit:

\begin{lemma}
\label{le:uniquefix}
Let $g$ be a leaf of $\cG^s$ or $\cG^u$ preserved by a non-trivial element $\gamma$ of $\pi_1(M)$. Then, $\gamma$ preserves one and only one element of $g$.
\end{lemma}

\begin{remark}
In the sequel, $\phi$ will always denote a codimension one Anosov action of $\RR^k$ on a closed manifold $M$ of dimensions $n+k$. All the statements we will prove are invariant through finite coverings and $\TT^\ell$-extensions. Therefore, we can always assume without loss of generality that the various foliations are oriented and transversely oriented, and that $\phi$ is irreducible.
\end{remark}

\section{Cross-section to one Anosov element}

\subsection{$1$-suspended Anosov actions}

\begin{definition}
\label{def:susf}
{\rm
Let $f: N \rightarrow N$ be a partially hyperbolic diffeomorphism whose central direction is tangent
to the orbits of a locally free action of $\RR^{k-1}$ commuting with $f$. The suspension of $f$ is the manifold $N_f$, quotient
of $N \times \RR$ by the diagonal action of $\ZZ$ defined by:
$$k.(x,t)=(f^k(x), t-k),$$
equipped with the action of $\RR^k \approx \RR^{k-1} \oplus \RR$ defined by:
$$(a,s).[x,t] = [a.x, t+s]$$
The manifold $N_f$ equipped with this action of $\RR^k$ is called a \textit{$1$-suspended Anosov action.}
}
\end{definition}

$1$-suspended Anosov actions are obviously Anosov actions of $\RR^k$, for which the element
of coordinates $(0,1)$ in $\RR^k \approx \RR^{k-1} \oplus \RR$ is Anosov. Observe that the submanifold $N \times \{ 0 \}$ projects
in $N_f$ to a cross section to the flow induced by the translation on the coordinate $t$; i.e. every orbit of this flow intersects
this submanifold infinitely many times.

Inversely:

\begin{proposition}
An Anosov action $(M,\phi)$ is topologically equivalent to a $1$-suspended Anosov action if and only if
it admits a $1$-cross-section, i.e. a closed codimension one submanifold $N \subset M$ which is a cross-section for the flow $\phi^{ta_0}$ generated by an Anosov element $a_0$.
\end{proposition}

\begin{proof}
Let $N$ be a cross-section to $\phi^{ta_0}$.
Let $f: N \rightarrow N$ be the first return map along the orbits of $\phi^{ta_0}$. This map is clearly partially hyperbolic,
admitting as central direction the intersection between the tangent bundle to the orbits of $\phi$ and the tangent bundle
of $N$.

Recall that ${\mathcal R}^k$ denotes the Lie algebra of $\RR^k$. Let $H$ be the image $\zeta({\mathcal R}^k)$, and let
$X_0$ be the element of $H$ generating the flow $\phi^{ta_0}$. For every $X$ in $H$ and every $x$ in $N$, let $\bar{X}(x)$ be the projection of $X(x)$ on $T_xN$ relatively to the direction defined by $X_0$. It defines a vector field
on $N$. Moreover, for a given $x$, the space formed by the tangent vectors $\bar{X}(x)$ is the intersection between $T_x\phi$ and $T_xN$.

There is a function $r_X: N \to \RR$ such that:
$$X(x) = \bar{X}(x) + r_X(x)X_0(x)$$
Now if $Y$ is another element of the Lie algebra $H$, we have:
$$0 = [X, Y] = [\bar{X}, \bar{Y}] + (\bar{X}(r_Y)-\bar{Y}(r_X))X_0$$
Hence, since $[\bar{X}, \bar{Y}]$ is tangent to $N$, it vanishes.

Therefore, the vector fields $\bar{X}$ commute altogether, and define a locally free action of $\RR^{k-1}$ on $N$.
Moreover, they all commute with the first return map $f$. The proposition is proved.
\end{proof}

\subsection{Schwartzman criterion}
\label{sub:schwartzman}
\begin{definition}
{\rm
Let $a$ be an element of $\RR^k$.  A {\it $a$-periodic orbit} is a loop $t \to \phi^{ta}(x)$ where $t$ describes $[0,1]$ and $x$ a fixed point of $\phi^a$.
}
\end{definition}

A {\it regular cone} is an open convex cone in $\RR^k$ containing only Anosov elements.
Every periodic orbit $c$ defines an element $[c]$ of $H_1(M, \RR)$, on which every element of $H^1(M, \RR)$ can be evaluated.

\begin{theorem}
\label{teo:schwarz}
Let $\cC$ be a regular cone in $\RR^k$.
Let $\epsilon$ be a positive real number and $\omega$ be an element of $H^1(M,\RR)$ such that, for every  $a$-periodic orbit $c$ with $a$ in $\cC$, one has
$\omega([c]) > \epsilon \Vert a \Vert_0$.
Then, $(M, \phi)$ is topologically equivalent to an $1$-suspended Anosov action of $\RR^k$.
\end{theorem}

\begin{proof}
We fix an element $a_0$ of $\cC$, and aim to prove that the flow $\phi^{ta_0}$ admits a cross-section. To that purpose,
we use the following criterion by Schwartzman: in \cite{schwartzman}, Schwartzman associates to every
$\phi^{ta_0}$-invariant borelian measure a homology class $A_\mu$ in $H^1(M, \RR)$, called the asymptotic cycle.

$A_\mu$ depends continuously and linearly on $\mu$. When $\mu$ is ergodic, $A_\mu$ is defined in the following way:
let $x$ be $\mu$-regular point, i.e. a point for which, for every continuous function $h: M \rightarrow \RR$ one has:
$$\int h d\mu = \lim_{T \to \infty} \frac1T \int^T_0 h(\phi^{ta_0}(x))dt$$
Then $A_\mu$ is the limit in $H_1(M, \RR)$ of $\frac1T [c_T]$, where $c_T$ is the trajectory  $t \to \phi^{ta}(x)$ closed
by a small path, of length less than the diameter of $M$, joining $x$ to $\phi^{Ta_0}(x)$. Indeed, since we divide by $T$, the limit
does not depend on the choice of this additional path. Moreover, since $\mu$ is ergodic, the limit does not depend
on the choice of the regular point.

In  \cite{schwartzman}, Schwartzman  proved that if there is some element $\omega \in H^1(M, \RR)$ for which all the pairings
$\omega(A_\mu)$ is positive, for every $\phi^{ta_0}$-invariant measure $\mu$, then $\phi^{ta_0}$ admits a cross-section.

Since the compact set of all invariant \textit{probability} measures is the closed convex hull of ergodic invariant probability measure,
it is enough for us to show that there exists $\epsilon' > 0$ such that $\omega(A_\mu) \geq \epsilon'$ for all ergodic probability measures $\mu$.

Now, in that situation, one can select as $\mu$-regular point an element $x$ which is also $\phi^{a_0}$-recurrent: there
is an increasing sequence $T_n$ converging to $+\infty$ such that $d(x, \phi^{T_na_0}(x)) \leq 1/n$.
By the closing Lemma (cf. \cite[Theorem 2.4]{katokihes}), there exist positive constants $C$, $\lambda$, a
sequence of points $y_n$ and differentiable maps $\gamma_n:[0, T_n]\to \mathbb{R}^k$
  such that for all $t\in [0, T_n]$ we have
  \begin{enumerate}
    \item $d(\phi^{ta_0}(x),\phi^{\gamma_n(t)}(y_n))<Ce^{-\lambda(\min\{t,T_n-t\})}d(\phi^{T_na_0}(x), x);$
    \item $\phi^{\gamma_n(T_n)}(y_n)=\phi^{\delta}(y_n)$ where $\|\delta\| < C d(\phi^{T_na_0}(x), x);$
    \item $\|\gamma_n'-a_0\|< Cd(\phi^{T_na_0}(x), x).$
  \end{enumerate}

In each of these inequalities, the right term is bounded from above by $C/n$. Observe that by item $(3)$, every $\gamma_n'(t)$, for
$n$ big enough, belong to $\cC$ for all $t$, and in particular, $\gamma_n(T_n)$ lies in $\cC$.

Let us closes the trajectory $t \to \phi^{ta}(x)$
by the small geodesic of length $\leq C/n$ joining $x$ to $\phi^{T_na_0}(x)$, and the path $t \to \phi^{\gamma_n(T_n)}(y_n)$
by a small path of length $\leq C/n$ contained in the $\phi$-orbit. According to item $(1)$, these two loops are homotopic;
thus the pairing $\omega([c_{T_n}])$ is equal to the evaluation of $\omega$ on the homology class defined by the second.
But since this second loop is contained in a $\RR^k$-orbit, it is homotopic to the $(\gamma_n(T_n)-\delta)$-periodic orbit
of $y_n$ (cf. item $(2)$). By hypothesis, the result is bounded from below by $\epsilon \Vert \gamma_n(T_n)-\delta \Vert_0$.

Therefore, $\omega(A_\mu)$, which is the limit of $\omega([c_T])/T$, is bounded from below by the limit of
$\epsilon\frac{\Vert \gamma_n(T_n)-\delta \Vert_0}{T_n}$. But according to item $(3)$:
$$\Vert \gamma_n(T_n) - T_na_0 \Vert_0 \leq T_nC/n$$
Hence:
\begin{eqnarray*}
\frac{\Vert \gamma_n(T_n)-\delta \Vert_0}{T_n} & \geq & \frac{\Vert \gamma_n(T_n)\Vert_0}{T_n} - \frac{\Vert \delta \Vert_0}{T_n}\\
                                                                                    & \geq & \frac{\Vert T_na_0 \Vert_0}{T_n} - \frac{\Vert \gamma_n(T_n)-T_na_0 \Vert_0}{T_n} - \frac{C}{nT_n}\\
                                                                                    & \geq & \Vert a_0 \Vert  - \frac{C}n - \frac{C}{nT_n}
\end{eqnarray*}
Hence, taking $\epsilon' = \epsilon\Vert a_0 \Vert_0/2$ one has, for $n$ big enough:
$$\omega(A_\mu) \geq \epsilon'$$
The Theorem is proved.
\end{proof}

\section{Anosov action with tranversely affine stable foliation}
In this section, we will prove a strong version of Theorem~\ref{teo:transaffine} stated in the introduction:

\begin{theorem}
\label{teo:transaffinecomplet}
Let $\phi$ be a codimension one Anosov action on a closed manifold $M$. Then the following properties are equivalent:
\begin{enumerate}
\item the fundamental group of $M$ is solvable,
\item the weak stable foliation $\cF^s$ admits a $C^0$ transverse affine structure,
\item the weak stable foliation $\cF^s$ admits a $C^1$ transverse affine structure,
\item $\phi$ is splitting.
\end{enumerate}
Moreover, if one of these hypothesis are satisfied, then the transverse affine structure is complete, i.e. the leaf space $\cL^s$
is affinely isomorphic to the affine line $\RR$, and the holonomy representation $\rho^s: \pi_1(M) \rightarrow \operatorname{Aff}(1, \RR)$ is faithfull.
\end{theorem}

Observe that according to Theorem \ref{thm:ReducingAction} the proof of Theorem \ref{teo:transaffine} can be reduced
to the case of irreducible actions\footnote{It seems that since main part of our arguments are valid up to 
finite coverings, one should replace in item $(1)$ <<solvable>> by
virtually solvable>>, but see Remark \ref{rk:pasvir}.} (concerning item (1), see Remark \ref{rk:solvablered}).
The proof is divided in various subsections.

\subsection{Completeness of transverse affine structures}
Let $(M, \phi)$ be an irreducible codimension one action whose weak stable foliation $\cF^s$ is $C^0$ transversely affine.
This transverse affine structure induces on each strong unstable leaf an affine structure.

\begin{lemma}
\label{le:ucomplete}
Each leaf of $\cF^{uu}$ is complete, i.e. affinely isomorphic to $\RR$.
\end{lemma}

\begin{proof}
For every $x$ in $M$, let $y(x)$ be the point $y$ in $\cF^{uu}(x)$ at distance $1$ from $x$ along the leaf and such
that $y > x$ for the orientation along the leaf. Then,
there is a unique affine map $f_x$ sending the leaf $\cF^{uu}(x)$ into $\RR$, such that $f_x(x)=0$ and $f_x(y(x))=1$.
The image of $f_x$ is an interval $]-a(x), b(x)[$ with $0 < a(x), b(x) \leq +\infty$. The functions $a: M \rightarrow \RR_+^\ast \cup \{ +\infty \}$ and
$b: M \rightarrow \RR_+^\ast \cup \{ +\infty \}$ are upper semi-continuous. Hence they both admit a positive lower bound. Assume that
one of them, let say $b$, has a finite lower bound $b_0$, and let $x$ be a point where $b(x) < b_0 + 1$. Consider now an element $h_0$
of the Anosov cone $\cA_0$. For every $t > 0$, $\phi^{-th_0}$ is an affine transformation between $\cF^{uu}(x)$ and $\cF^{uu}(\phi^{-th_0}(x))$.
It follows that $f_{\phi^{-th_0}(x)}  \circ \phi^{-th_0} = \mu_tf_x$, for some $\mu_t > 0$. Now, since $\phi^{-h_0}$ contracts exponentially strong stable
leaves, the distance between $\phi^{-th_0}(x)$ and $\phi^{-th_0}(y(x))$ decreases exponentially with $t$. In other words,
there are $C, \lambda > 0$ such that $f_{\phi^{-th_0}(x)}  \circ \phi^{-th_0}(y(x)) \leq C\exp(-\lambda t)$. It follows that $\mu_t \leq C\exp(-\lambda t)$.
Hence, $b(\phi^{-th_0}(x)) = \mu_t b(x) \leq C\exp(-\lambda t)(b_0+1)$. This is impossible since $b$ admits a positive lower bound $b_0$.

This contradiction shows that $a$ and $b$ have only infinite value, i.e. that the leaves of $\cG^u$, as affine manifolds, are complete affine lines.
\end{proof}

\begin{proposition}
\label{pro:split}
The Anosov action is splitting.
\end{proposition}

\begin{proof}
Let $\cD: \widetilde{M} \rightarrow \RR$ be the developping map of the transverse affine structure (i.e. a map constant along the leaves
of $\widetilde{\cF}^s$ and which is affine with respect to the transverse affine structure and the usual affine structure on $\RR$).
Let $F_0$ be a leaf of $\widetilde{\cF}^s$, and let $U_0$ be the $\widetilde{\cF}^{uu}$-saturation of $F_0$. Let $F_1$ be a leaf of $\widetilde{\cF}^s$
containing a point of $U_0$.

Let $E_1$ be the set of leaves of $\widetilde{\cF}^s$ whose image by $\cD$ is equal to $\cD(F_1)$. For every $F$ in $E_1$, let $V(F)$ be the open
subset  of $F_0$ made of points whose unstable leaf intersects $F$.

According to Lemma \ref{le:ucomplete}:

$$F_0 = \bigcup_{F \in E_1} V(F)$$

Moreover, the restriction of $\cD$ to any unstable leaf is injective. It follows that if $F$, $F'$ are two different
elements of $E_1$, then $V(F)$ and $V(F')$ are disjoint.

Since $F_0$ is connected, it follows that $F_0 = V(F_1)$. In other words, $U_0$ is $\widetilde{\cF}^s$-saturated.
Since it is by definition $\widetilde{\cF}^{uu}$-saturated, $U_0$ is the entire universal covering $\widetilde{M}$.
The proposition is proved.
\end{proof}

\begin{corollary}
\label{cor:detailaff}
The transverse affine structure is complete, i.e. the developping map $\cD: \widetilde{M} \rightarrow \RR$
is a local  fibration whose fibers are exactly the leaves of $\widetilde{\cF}^s$. In particular,
$\cL^s$ is diffeomorphic to $\RR$, and the holonomy morphism $\rho^s: \pi_1(M) \rightarrow \operatorname{Aff}(1, \RR)$ is injective.
\end{corollary}

\subsection{Smoothness of transverse affine structure}

The foliation $\cF^s$ is $C^1$, hence it makes sense for a transverse affine structure to be $C^1$, but not to be $C^2$.
In \cite{Bar-Maq} we proved that among $C^2$ affine structures along the leaves of
$\cF^{uu}$ varying continuously with the leaf, there is one and only one that is preserved by $\phi$.

It is not obvious that this affine structure coincide everywhere
with the one induced by the affine structure transverse to $\cF^s$, since the uniqueness
result proved in \cite{Bar-Maq} works only among affine structures along $\cF^{uu}$ at least $C^2$, whereas the
transverse affine structure constructed in Corollary \ref{cor:transaffine} is \textit{a priori} only $C^0$.

\begin{proposition}
\label{pro:C0C1}
Let $\phi$ be a codimension one Anosov action on a manifold $M$. Then, any $C^0$ affine structure transverse to
$\cF^s$ is actually $C^1$.
\end{proposition}

\begin{proof}
We have to prove that the developping map $\cD: \widetilde{M} \rightarrow \RR$ is $C^1$.
Let $x$ be a point in $\widetilde{M}$ whose projection in $M$ has a compact $\phi$-orbit. Then, according to Lemma \ref{le:ucomplete} and Proposition \ref{le:ucomplete},
the restriction of $\cD$ to $\widetilde{\cF}^{uu}(x)$ is onto. Since the leaf $\widetilde{\cF}^{uu}(x)$ is smooth (even if the foliation itself is only
h\"{o}lder), it is enough to prove that this restriction is $C^1$.

The projection by $\pi^\phi$ of $\widetilde{\cF}^{uu}(x)$ in $Q^\phi$ is a leaf $g^u$ of $\cG^u$. There is a map $d: g^u \rightarrow \RR$
such that $\cD = d \circ \pi^\phi$ along $\widetilde{\cF}^{uu}(x)$. Since $\pi^\phi$ is smooth, we have to show that $d$ is $C^1$.

Assume $k\geq 2$.
Let $H$ be the fundamental group of the $\phi$-orbit of $x$, i.e. the stabilizer of $\pi^\phi(x)$ of the $\Gamma$-action.
Then $H$ is a free abelian group, isomorphic to $\ZZ^k$, preserving $g^u$. Moreover, the $C^2$-affine structure along $\widetilde{\cF}^{uu}(x)$
constructed in \cite{Bar-Maq} provides a $C^2$ parametrization of $g^u$ by $\RR$, so that for this parametrization,
the action of $H$ is an action by scalar multiplications, through a morphism $\rho^c: H \rightarrow \RR_\ast^+$:
$$h(t) = \rho^c(h) . t$$
where $.$ denotes the usual multiplication in $\RR$.

Now since $d$ comes from the developping map, we have for every $t$ in $\RR$:
$$d(\rho^c(h).t) = \rho^s(h).d(t)$$
where $\rho^s$ is the holonomy of the affine structure transverse to $\cF^s$ (here, we choose the affine parametrization of the target of $d$
so that $d(\pi^\phi(x)) = 0$; the affine action of $\rho^s(H)$ thus preserves $0$ and is indeed a linear action by scalar multiplications).

Consider the map $\log \circ d \circ \exp: \RR \rightarrow \RR$: it is a topological conjugacy between the $H$-actions by translations on $\RR$, one
given by $\log \circ \rho^c$, the other by $\log \circ \rho^s$. In other words, it is an order preserving map between $\log \circ \rho^c(H)$ and
$\log \circ \rho^s(H)$. Since here we assume $k \geq 2$, these two groups are dense subgroups of $\RR$. It follows that $\log \circ d \circ \exp$
is actually the multiplication by a constant positive scalar, and thus, that the restriction of $d$ to $]0, +\infty[$ has the form:
$$d(t)=t^\alpha, \;\; \operatorname{  for } t > 0$$
Similarly, there is a positive real number $\beta > 0$ such that:
$$d(t)=-(-t)^\beta, \;\; \operatorname{  for } t < 0$$

Hence it proves that $d$ is $C^1$, at least out of $0$. It follows that $\cD$ is $C^1$, except maybe along $\widetilde{\cF}^s(x)$.
To prove that $\cD$ is $C^1$ also near $\widetilde{\cF}^s(x)$ is just a matter to repeat the argument for another point $x'$ projecting
in $M$ with compact orbit but with $\widetilde{\cF}^s(x') \neq \widetilde{\cF}^s(x)$. It finishes the case $k\geq 2$.

The case $k=1$ is similar: in that case, according to \cite{plante2, matsumotosolv}, we know that $\phi$ admits a global section $N$, and that the first return
map is an Anosov diffeomorphism of codimension one. Hence, according to \cite[Corollary (1.4)]{newhouse} the global section is homeomorphic to the torus $\TT^n$,
and the first return map is topologically conjugate to a linear hyperbolic automorphism. The problem is then reduced
to prove that the affine structure transverse to the foliation induced on $N$ by $\cF^s$ is $C^1$. The proof is similar to the proof in the previous case
$k\geq2$, but by considering the (free) action of $\pi_1(N) \approx \ZZ^n$ on the leaf space instead of the action of the stabilizer of a leaf.
\end{proof}

\subsection{Proof of Theorem \ref{teo:transaffinecomplet}}
The implication $(1) \Longrightarrow (2)$ is Corollary 2 in \cite{matsumotosolv}. $(2) \Longrightarrow (3)$ is the content of
Proposition \ref{pro:C0C1}. $(3) \Longrightarrow (4)$ is Proposition \ref{pro:split}. The final implication $(4) \Longrightarrow (1)$
is due to the following Theorem of V.V. Solodov (cf. \cite[Theorem 6.12]{ghysolodov}):

\begin{theorem}
\label{thm:solodov}
Let $\Gamma$ be a non-abelian group of orientation preserving homeomorphisms of the line such that
every element has at most one fixed point. Then, $\Gamma$ is isomorphic to a subgroup
of the affine group $\operatorname{Aff}(1, \RR)$ of the real line, and the action of $\Gamma$ on the line is semi-conjugate to the corresponding affine action.
\end{theorem}

Indeed, we can apply this theorem to the group $\Gamma = \pi_1(M)$ acting on the leaf space $\cL^s \approx \RR$: assume that some element $\gamma$ of $\pi_1(M)$
admits two different fixed points in $\cL^s$; i.e. preserves two leaves $\cG^s(x)$, $\cG^s(x'),$ where $x$, $x'$ are $\gamma$-fixed points in $Q^\phi$. Then, since the action is splitting,
the unstable line $\cG^u(x)$ intersects these two leaves: the first at $x$, the other at some point $y$. Now, since $x$ is fixed by $\gamma$, the leaf $\cG^u(x)$ is preserved
by $\gamma$. Therefore, $y$ is also a $\gamma$ fixed point. It implies that $\cG^u(x)$ contains two $\gamma$-fixed points, which is impossible (cf. Lemma~\ref{le:uniquefix}). Hence, according
to Theorem \ref{thm:solodov}, $\pi_1(M)$ is solvable.

The additional statement in Theorem \ref{teo:transaffine}  about the completeness of transverse affine structures is Corollary \ref{cor:detailaff}.

\begin{remark}
\label{rk:pasvir}
Up to that point, since our proof is up to finite coverings, we have actually shown Theorem \ref{teo:transaffinecomplet} where item $(1)$ is 
replaced by item $(1')$: \textit{$\pi_1(M)$ is virtually solvable,}. However, as a matter of fact, $\pi_1(M)$ is solvable as soon as it is virtually solvable. Indeed: assume that $(1')$ is true. 
Applying the result a finite covering of $M$, we get that the action is splitting. Let now $\Gamma_0$ be the index $2$ subgroup
of $\pi_1(M)$ corresponding to elements preserving the orientation of $\cL^s$, ie. to the double covering where the weak stable foliation is transversely orientable.
Then Solodov's Theorem implies that $\Gamma_0$ is solvable. Since the commutator subgroup of $\pi_1(M)$ is obviously contained in $\Gamma_0$, it
follows that $\pi_1(M)$ is solvable too.
\end{remark}

\subsection{Splitting Anosov actions are $1$-suspensions}

\begin{theorem}
\label{teo:splitsusp}
Let $\phi$ be a codimension one Anosov action of $\RR^k$. Assume that $\cF^s$ is transversely affine.
Then, $\phi$ is topologically equivalent to a $1$-suspended Anosov action.
\end{theorem}

\begin{proof}
One defines, exactly as in \cite[\S 2]{plante2}, a cohomology class $\omega$ expressing the dilatation of the transverse affine structure along stable leaves.
More precisely: since $\cF^s$ is transversely affine, there is a covering $\cU$ of $M$ by $\cF^s$-distinguished subsets $U_i$: for each $i$, the leaves of
the restriction of $\cF^s$ to $U_i$ are the fibers of a submersion $\pi_i: U_i \rightarrow \RR$; and the transition functions $\tau_{ij} = \pi_j \circ \pi_i^{-1}$ are affine maps.
Choose the covering $\cU$ sufficiently fine so that $\pi_1(M) \approx \pi_1(N(\cU))$ where $N(\cU)$ is the
nerve of the covering $\cU$ and let $\gamma$ be an element of $\pi_1(M)$ which is represented by a
chain of open sets $(U_1,...,U_m)$ in $\cU$ ($U_i \cap U_{i+1} \neq \emptyset$ for  $i = I, ...,n-1$ and $U_n \cap U_1 \neq \emptyset$).
Denote by $J_i$ the Jacobian of the transition function $\pi_{i+1} \circ \pi_i^{-1}$ for
$i = 1,...,n -1$ and by $J_n$ the Jacobian of $\pi_{1} \circ \pi_n^{-1}$. The $J_i$ are constants since $\cF^s$ is
transversely affine and we define $\omega(\gamma) = \sum_{i=1}^n \log \|J_i \|$. It follows easily from the cocycle
property of transition functions that $\omega(\gamma)$ depends only on the homotopy class of $\gamma$ so
that $\omega$ is a well-defined element of $\operatorname{Hom}(\pi_1(M), \RR) \approx H^1(M, \RR)$.

Let now $a_0$ be an element Anosov of $\RR^k$.
Then, there is a constant $\lambda > 1$ and a small neighborhood $U$ of $a_0$ in $\RR^k$ such that, for every
$a$ in $U$, for every $\phi^a$ fixed point $x$ and every $u$ in $E^{uu}(x)$, we have $\Vert D\phi^a(u)\Vert \geq  \lambda \Vert u \Vert$. This implies that the
Jacobian of the restriction of $\phi^a$ to the strong unstable manifold at $x$ is $\geq \lambda$.
It follows, reducing $\lambda$ if necessary, that there is a regular convex cone $\cC$ containing $a_0$ such that, for every $a$ in $\cC$ the Jacobian of
the restriction of $\phi^a$ to the strong unstable manifold at $x$ is $\geq \lambda^{\Vert a \Vert_0}$.

Since
this Jacobian is unchanged by a differentiable change of coordinates we conclude
that $\omega([c]) \geq \log\lambda\Vert a \Vert_a$ for any periodic orbit $c$ of an element $a$ of $\cC$. This concludes the proof
of Theorem \ref{teo:splitsusp} thanks to Theorem \ref{teo:schwarz}.
\end{proof}

\begin{remark}
In this proof, we used the fact that the transverse affine structure is not only $C^0$, but $C^1$ (cf. Proposition \ref{pro:C0C1}).
See the comment in \cite{plante2} just after the statement of Theorem B.
\end{remark}

\section{Actions transverse to a fibration}
\label{sec:fibration}
In this section, we assume that $\phi$ is transverse to a locally trivial fibration $\pi: M \rightarrow B$. Let us start with an easy observation:

\begin{lemma}
\label{le:basetorus}
Up to a finite covering, one can assume that the basis $B$ is diffeomorphic to $\TT^k$.
\end{lemma}

\begin{proof}
Let $\cO_0$ be a compact orbit of the action. It is transverse to the fibers of $\pi$; hence
the restriction of $\pi$ to $\cO_0$ is a covering map $\pi_0$. Let $\pi'$  be the
pull-back $\pi$ by $\pi_0$, i.e. the unique map $\pi': M' \rightarrow \cO_0$ such that the following
diagram commutes:

$$ \xymatrix{
    M' \ar[r]^{\pi'} \ar[d]_{\pi'_0} & \cO_0 \ar[d]^{\pi_0} \\
    M \ar[r]^\pi & B
  }$$

Then $\pi'$ is a locally trivial fibration, the map $\pi'_0: M' \rightarrow M$ is a finite covering, and
the action of $\RR^k$ on $M$ lifts to an Anosov action of $\RR^k$ on $M'$ transverse to $\pi'$. The lemma
follows since $\cO_0$ is diffeomorphic to $\TT^k$.
\end{proof}

From now, we assume that the basis $B$ is the torus $\TT^k$.

\begin{corollary}
\label{cor:susmonodromy}
The Anosov action $(M, \phi)$ is topologically equivalent to the suspension of an action of $\ZZ^k$ on a manifold $F_0$.
\end{corollary}

\begin{proof}
Let $F_0$ be a fiber of $\pi$.
The orbits of $\phi$ are the leaves of a foliation $\cF_\phi$  on $M$ transverse to the fibers of $\pi$. Hence, $(M,  \cF_\phi, \pi)$
is a \textit{foliated bundle}, or a \textit{flat bundle}. In other words, there is a diffeomorphism between $(M,  \cF_\phi, \pi)$
and the suspension $(M, \cF_\rho, \pi_\rho)$ of a representation $\rho: \pi_1(\TT^k)\approx \ZZ^k  \rightarrow \operatorname{Diff}(F_0)$.

The action of $\RR^k$ on $\TT^k$ by translations lifts in a unique way to a locally free action of $\RR^k$ preserving each orbit of $\cO$.
This new action is obviously topologically equivalent to $\phi$ and is also isomorphic to the suspension
of $\rho: \ZZ^k \rightarrow  \operatorname{Diff}(F_0)$.
\end{proof}

\begin{remark}
\label{rk.basepoint}
In particular, there is a global $\rho(\ZZ^k)$-fixed point, which is the intersection between the fiber $F_0$ and the compact orbit $\cO_0$.
In the sequel, we use it as a base point in $F_0$, and denote it by $x_0$. We select a lifting $\tilde{x}_0$ of $x_0$ in the universal covering
$\widetilde{F}_0$ of $F_0$. Then, we can lift $\rho$ to a morphism $\tilde{\rho}: \widetilde{F}_0 \rightarrow \widetilde{F}_0$, uniquely characterized
by the requirement that every $\tilde{\rho}(\gamma)$ preserves $\tilde{x}_0$.
\end{remark}

\subsection{The monodromy cocycle}
The key difficulty is that the initial Anosov action $\phi$ and the new action $\phi_0$, which is a suspension, may not coincide.
We can define a (smooth) cocyle $\eta: M \times \RR^k \rightarrow \RR^k$ in the following way: $\eta(x, u)$ is the unique element $a$ of
$\RR^k$ such that $\phi^u_0(x)=\phi^{a}(x)$ (the ambiguity arising when the orbit of $x$ is periodic is surrounded by defining $\eta$
at the universal covering level). Then, for every $\gamma$ in $\ZZ^k \subset \RR^k$, the monodromy $\rho(\gamma)$ is the map from $F_0$ onto
itself mapping $x$ on $\phi^{\eta(x, \gamma)}(x)$.

\begin{definition}
{\rm
An element $\gamma$ of $\ZZ^k$ is \textit{Anosov} if for every $x$ in $F_0$ the cocycle $\eta(x, \gamma)$ belongs to the codimension one Anosov cone $\cA_0$.
}
\end{definition}

In the case $k=1$, the existence of Anosov elements in $\ZZ^k$ is obvious since in that case $\cA_0$ is the half-line $\RR^+_\ast$. It is well-known
that in that case, the action $\phi$ is topologically (even H\"older) conjugate to the suspension of a toral automorphism.

But the case $k \geq 2$ is far from obvious.
As a matter of fact, we are unable to prove in a whole generality that the new action $\phi_0$ is Anosov, even less that
$\phi_0$ is a codimension one Anosov action! However:

 \begin{lemma}
 \label{le:anosovmonodromy}
 If $k\geq2$ and $\ZZ^k$ contains an Anosov element, then $\phi$ is H\"older conjugate to the suspension of a $\ZZ^k$-action by toral automorphisms up to an automorphism of $\RR^k$.
 \end{lemma}

 \begin{proof}
 Let $\gamma$ be an Anosov element of $\ZZ^k$. Then, since for every $x$ in $F_0$, $\eta(x,\gamma)$ is Anosov, it follows easily that the monodromy $\rho(\gamma)$
 is an Anosov diffeomorphism of $F_0$, whose unstable direction has dimension one. According to Franks-Newhouse Theorem, it is topologically conjugate
 to a hyperbolic linear transformation of the torus, where the topological conjugacy is H\"older continuous.  Now it is well-known that any homeomorphism commuting
 with a linear hyperbolic automorphism of the torus is itself
 a linear automorphism (see for example \cite[Lemma 2.4]{KL2}; observe that $x_0$ is a common $\rho(\ZZ^k)$-fixed point). In other words, $\phi$ is H\"older topologically equivalent to the suspension of a $\ZZ^k$-action by toral automorphisms.
 According to \cite[Theorem 2.12 (b)]{katokihes}, $\phi$ is then H\"older conjugate to this suspension up to an automorphism of $\RR^k$.
 \end{proof}


 \vspace{.5cm}

\begin{remark}
Remind that we can expect an even stronger conclusion: the conjugacy might be smooth, not only H\"older (\cite[Conjecture 1.1]{kalininspatz}).
\end{remark}

\subsection{Anosov actions transverse to a fibration are splitting}
We don't know how to prove, for a general codimension one Anosov action of $\RR^k$ transverse to a fibration, the existence of an Anosov monodromy element.
Nevertheless, we can establish many properties in this general case; and the remaining part of this section is devoted to the proof of these properties.

Since $F_0$ is transverse to $\phi$, it is transverse to the weak foliations $\cF^s$, $\cF^u$.
We denote by $\hF^s$, $\hF^u$ the foliations induced on $F_0$ by $\cF^s$, $\cF^u$.
Every monodromy transformation $\rho(\gamma)$ for $\gamma$ in $\ZZ^k$ preserves the foliations $\hF^s$, $\hF^u$. Moreover:

\begin{lemma}
\label{le:propseq}
Let $\widetilde{F}_0$ be the universal covering of $F_0$, and let $\Gamma_F$ be the fundamental group of $F_0$.
Then the inclusion $F_0 \subset M$ induces an inclusion $\Gamma_F \subset \Gamma$, and lifts to an inclusion $\widetilde{F}_0 \subset \widetilde{M}$.
The restriction of $\pi^\phi:  \widetilde{M} \rightarrow Q^\phi$ to $\widetilde{F}_0$ is a $\Gamma_F$-equivariant diffeomorphism,
mapping the liftings of $\hF^s$, $\hF^u$ on the foliations $\cG^s$, $\cG^u$.
\end{lemma}

\begin{proof}
Since $B$ is diffeomorphic to $\TT^k$, the second homotopy group $\pi_2(B)$ is trivial. The fibration  exact sequence for $\pi$ provides the following exact sequence:
$$0 \rightarrow \pi_1(F_0) \rightarrow \pi_1(M) \rightarrow \ZZ^k \rightarrow  0 \;\;\; (\ast)$$
Hence we get the inclusion $\Gamma_F \subset \pi_1(M)$. Moreover, it follows that every connected
component of $p^{-1}(F_0)$ in $\widetilde{M}$ is homeomorphic to $\widetilde{F}_0$.

By construction, $M_\rho$ is the quotient of $F_0 \times \RR^k$ by the action of $\ZZ^k$ defined by:
$$h.(x, v) = (\rho(h)x, v+h); \;\; h \in \ZZ^k$$
Hence, the universal covering of $M \approx M_\rho$ is the product $\widetilde{F}_0 \times \RR^k$,
where the orbits of the lift of $\phi$ have the form $\{ \ast \} \times \RR^k$. The lemma follows.
\end{proof}

\begin{lemma}
\label{le:R}
Every leaf of $\hF^s$ (respectively $\hF^u$) is diffeomorphic to $\RR^{n-1}$ (respectively $\RR$).
\end{lemma}

\begin{proof}
Assume that some leaf $\hat{L}$ of $\hF^s$ contains a loop $c$ homotopically non-trivial.
According to Lemma \ref{le:propseq}, $c$ is also homotopically non-trivial in $M$.
Observe that $c$ is a loop inside the leaf of $\cF^s$ containing $\hat{L}$.
Since $\phi$ is irreducible, it follows that this leaf is the stable leaf of a compact orbit $\cO_p$
and that $c$ is homotopic to a loop in $\cO_p$. At the one hand, $\pi_\ast([c])$ vanishes, since
$c$ is contained in a fiber, but on the other hand, it is non-zero since the restriction of $\pi$ to
$\cO_p$ is a covering map. Contradiction.

Therefore, every leaf of $\hF^s$ is simply connected, and thus is diffeomorphic to any lift of it
in $\widetilde{M}$. According to Lemma \ref{le:propseq} it is diffeomorphic to a leaf
of $\cG^s$, hence to $\RR^{n-1}$ (cf. \cite{Bar-Maq}).

The proof for $\hF^u$ is similar.
\end{proof}

\begin{corollary}
\label{cor:transaffine}
The foliation $\cF^s$ is transversely affine.
\end{corollary}

\begin{proof}
According to Lemma \ref{le:R} $\hF^s$ has no holonomy. It is well-known Theorem
by Sacksteder that in this case the
leaf space $\cL^s$ of the lift of $\hF^s$ in $\widetilde{F}_0$ is homeomorphic to the real line, and that
the fundamental group $\pi_1(F_0)$ is free abelian. Hence $\Gamma$ is solvable (even metabelian).
The corollary follows.
\end{proof}

\begin{proposition}
\label{pro:homotopytorus}
The fiber $F_0$ is homotopically equivalent to the torus $\TT^n$. If $n \neq 4$, $F_0$ is homeomorphic to $\TT^n$.
\end{proposition}

\begin{proof}
According to Theorem \ref{teo:transaffine}, we have a faithfull holonomy representation $\rho^s: \pi_1(M) \rightarrow \operatorname{Aff}(1, \RR)$.
It follows from Lemmas \ref{le:propseq} and \ref{le:R} that $\rho^s(\gamma)$ for $\gamma$ in $\Gamma_F = \pi_1(F_0)$ is a translation
of the real line. Recall the exact homotopy sequence:
$$0 \rightarrow \pi_1(F_0) \rightarrow \pi_1(M) \rightarrow \ZZ^k \rightarrow  0 \;\;\; (\ast)$$

Hence, any element $\gamma$ in $\pi_1(M)$ which is not in $\Gamma_F$ projects to a non-trivial element of $\ZZ^k$, i.e. representing a homotopically non-trivial loop
in $\cO_0$. The stable leaf containing $\cO_0$ has a non-trivial holonomy; it follows that $\rho^s(\gamma)$ is a pure homothety, not a translation.

Hence, $\Gamma_F$ is precisely the subgroup of $\Gamma$ made of elements $\gamma$ such that $\rho^s(\gamma)$ is a translation.
In particular, $\Gamma_F$ is a free abelian group. It is also the fundamental group of a closed $n$-dimensional manifold $F_0$
whose universal covering is diffeomorphic to $\RR^n$ (according to Lemma \ref{le:R} and Palmeira's Theorem). Therefore, $\Gamma_F$
is isomorphic to $\ZZ^n$, and $F_0$ is a $K(\ZZ^n, 1)$, i.e., homotopically equivalent to $\TT^n$. The final assertion when $n \neq 4$
follows by a Theorem of C.T.C. Wall (\cite{wall}).
\end{proof}

 \section{Integrable Anosov actions}
 The main objective of this section is to introduce the notion of integrability for Anosov actions and prove the two results stated below.

 \begin{theorem}
 \label{thm:integrable}
  If the action $\phi$ is integrable, then it is topologically equivalent to a suspension of an Anosov action of $\ZZ^k$ by linear automorphisms on the tori $\TT^n.$
  \end{theorem}

  The second result provides a sufficient condition for integrability  of codimension one Anosov actions.

   \begin{theorem}
   \label{thm:cond-integrable}
   If the action $\phi$ is volume preserving and the subbundle $E^{ss}\oplus E^{uu}$ is of class $C^1$, then the action $\phi$ is integrable.
   \end{theorem}

   Recall that $\mathcal{R}^k$ denotes the Lie algebra of $\RR^k$ and $H$ is the image $\zeta(\mathcal{R}^k)$. For each $i=1,\dots,k,$ let $X_i$ be the vector field in $H$  generating the flow
 $\phi_i^{ta_i}$, where $a_1, \dots, a_k$ is a base
 of $\RR^k$ which is contained in the Anosov cone $\mathcal{A}_0.$

  \subsection{Integrable Anosov actions are suspensions}
 Here we prove the Theorem \ref{thm:integrable} . We say that the action $\phi$ is {\it integrable}, if the subbundle $E^{ss}\oplus E^{uu}$ is integrable.

  \subsubsection*{Foliations defined by a system of $1$-forms}
 An $k$-tuple $\Omega=(\omega_1,\dots,\omega_k)$ of $1$-forms is
  a \textit{system of rank} $k$ on $M$ if the map $\Omega:TM\to \mathbb{R}^k$
  has rank $k$  when restricted to any fibre of $TM$.
  The $\ker \Omega$ is a subbundle of $TM$ which is called the \textit{kernel of the system}
  $\Omega$. If the bundle $\ker \Omega$ is integrable, we say that $\Omega$ is a \textit{integrable
  system}.

  The following proposition is well known, see for example
  \cite{hector}.

  \begin{proposition}
  \label{prop:hector}
 The kernel of a system $\Omega=(\omega_1,\dots,\omega_k)$ of rank
  $k$ is integrable  if and only if
  $$
  d\omega_i\wedge \omega_1\wedge\dots\wedge\omega_k=0
  $$
  for each $i$ in $\{1,\dots,k\}$.
  \end{proposition}

 \begin{proposition}
 \label{prop:int-fibration}
 If $\phi$ is an Anosov action of $\RR^k$ on $M$ for which $E^{ss}\oplus E^{uu}$
 is integrable, then  $\phi$ is transverse to a locally trivial fibration $\pi: M \rightarrow \TT^k$ .
 \end{proposition}

 This proposition is a generalization of a result (Theorem 3.1
 \cite{plante1}) obtained by Plante for Anosov flows and for its proof we will use
 the following result.

 \begin{theorem}[Plante \cite{plante1}]
 \label{thm:plante}
 Let $\cF$ be a $C^1$ foliation of codimension one on $M$ and
 $\psi^t$ a $C^1$ flow on $M$ transverse to $\cF$ which preserves
 the leaf of $\cF.$ Then $\cF$ is determined by a closed 1-form.
 \end{theorem}

 \vspace{.5cm}
 \begin{proof}[Proof of Proposition \ref{prop:int-fibration}]
 Let $\cF_i$ be the codimension one $C^1$ foliation on $M$ tangent to 
 $E^{ss}\oplus E^{uu}\oplus \RR X_1 \oplus\cdots \oplus \RR X_{i-1}\oplus \RR X_{i+1}\cdots \RR X_{k}$. Since the flow
 $\phi_i^t$ is transverse to $\cF_i$ and preserve the leaves of $\cF_i$, by
 Theorem \ref{thm:plante},  the foliation $\cF_i$ is defined by a closed 1-form
 $\omega_i$. We can assume that $\omega_i(X_i)=1$ for all $i=1,\dots,k.$ The system $\Omega=(\omega_1,\dots,\omega_k)$ of closed one forms
 is of rank $k$ and defines a foliation of codimension $k$ which is transverse to
 action $\phi.$ By a classical argument by Tischler (see for example \cite[Theorem 2.10]{plante1}), each
 one form $\omega_i$ may be a $C^0$ approximated by a closed one form $\omega_i'$ having rational periods.
 Therefore, the system $\Omega'=(\omega_1',\dots,\omega_k')$ also defines  a codimension $k$
 foliation that is transverse to the orbits of $\phi$ and whose leaves are  compacts.
 We can assume that each one form $\omega_i'$  has integer periods since any non-zero real multiple of $\Omega'$  determines the same foliation.

 Let $x_0$ be a basepoint in $M$. The application $\pi:M\to \TT^k$ given by
 $$
 \pi(x)=\Big(\int_{x_0}^x \omega_1'({\rm mod}1),\dots, \int_{x_0}^x \omega_k'({\rm mod}1)\Big)
 $$
 defines a locally trivial fibration (the definition of $\pi$ is independent of path). Moreover, $M$ is a bundle  over $\TT^k$ whose fibers are leaves of the foliation defined by $\Omega'$. This finishes the proof.
 \end{proof}

 \begin{proof}[Proof of Theorem \ref{thm:integrable}]
 With the notations above, we fix a fiber  $F_0=\pi^{-1}(0)$, note that $x_0\in F_0$ because
 $\pi(x_0)=0$. Since the action $\phi$ is topologically equivalent to the suspension of the monodromy action
 $\phi_0:\pi_1(\TT^k, 0)=\ZZ^k \to {\rm Diff}^1(F_0)$, it is sufficient to prove that some $\gamma \in \pi_1(\TT^k,0)=\ZZ^k$
 acts on the fiber $F_0$ as an Anosov diffeomorphism.

 Without less of generality we can assume that the orbit by $x_0$ is compact having its isotropy group generated by the set $\{a_1,\dots,a_k\}$ and, in the proof of Proposition \ref{prop:int-fibration}, we can assume that $\omega_i'(X_i)=1$ and $\omega_i'(X_j)=0$  for every $i$ and $j$ in $\{1,\dots,k\}$ with $j\ne i.$
Now we claim that  $\gamma$, the element of $\pi_1(\TT^k, 0)=\ZZ^k $ which is represented by the loop $\{(t,\dots, t); \  \ t\in \RR/\ZZ\}$ , acts on $F_0$ as an Anosov diffeomorphism.  To simplify the notation we write $\gamma =\{(t,\dots, t); \  \ t\in \RR/\ZZ\}. $ In fact, let  $p\in F_0$ and
 $\widetilde{\gamma}:[0,1]\to M$ defined by $\widetilde{\gamma}(t)=\phi^{t(a_1+\cdots+a_k)}(p)$,
 then since
 $$
 \int_{p}^{\widetilde{\gamma}(p)}\omega'_i=\int_{0}^{t}\omega'_i(\widetilde{\gamma}(p))\cdot (X_1+\cdots X_k)(\widetilde{\gamma}(p))=\int_{0}^{t}\omega'_i(\widetilde{\gamma}(p))\cdot X_i(\widetilde{\gamma}(p))=t,
 $$
 we obtain that

 $$
\begin{array}{lll}
\pi(\widetilde{\gamma}(t))  & =  &   \Big(\int_{x_0}^p\omega'_1+\int_{p}^{\phi^{t(a_1+\cdots+a_k)}(p)}\omega'_1,\dots, \int_{x_0}^p\omega'_k+\int_{p}^{\phi^{t(a_1+\cdots+a_k)}(p)}\omega'_k\Big )\\

  &  = &  \Big(\int_{p}^{\widetilde{\gamma}(p)}\omega'_1,\dots, \int_{p}^{\widetilde{\gamma}(p)}\omega'_k\Big ) \\
  &  =  &\gamma(t).
 \end{array}
$$
This means that $\tilde{\gamma}$ is a lift of  $\gamma$, consequently we have that  $\phi_0(\gamma)(p)=\phi^{a_1+\cdots +a_k}(p)$. Moreover, by convexity of $\mathcal{A}_0$,  the element $a_1+\cdots +a_k\in \mathcal{A}_0, $ i.e. is Anosov. Therefore
$\phi_0(\gamma)$ is an Anosov diffeomorphism and the proof is concluded.
\end{proof}

\subsection{A condition for integrability of codimension one Anosov actions}

This subsection is devoted to prove the Theorem \ref{thm:cond-integrable}, this will be done using some arguments  of Ghys in \cite{ghys}.
Recall that we are considering a codimension one Anosov  action of $\RR^k$ on a manifold $M$
of dimension at least $k+3$ which is volume preserving and such that the subbundle $E^{ss}\oplus E^{uu}$ is of class $C^1$.

\begin{proof}[Proof of Theorem \ref{thm:cond-integrable}]
For all $i\in \{1,\dots,k\}$ we consider the differential one form $\omega_i$ which is equal to $0$ on $E^{ss}\oplus E^{uu}\oplus \RR X_1 \oplus\cdots \oplus \RR X_{i-1}\oplus \RR X_{i+1}\cdots \RR X_{k}$ and to $1$ on the vector field $X_i$.  Note that the kernel of the system $\Omega=(\omega_1,\dots,\omega_k)$ is equal to $E^{ss}\oplus E^{uu}$.  It is clear that $\omega_i$ and
$d\omega_i$ are invariants by $D\phi^{ta_i}$ and consequently, by commutativity are also invariants
by $D\phi^{ta_j}$ for all $j\in \{1,\dots,k\}.$

We will show that $d\omega_i=0$ for all $i$, this will conclude the proof of the theorem.
We consider the 2-form $d\omega_1,$ the other cases are similar.

{\it Claim 1: the subbundle $T\phi$ is contained in the kernel of $d\omega_1$.}
Note that $T\phi$ is equal to $\RR X_1 \oplus\cdots \oplus \RR X_{k}$, hence it is sufficient to prove
that all the vector fields $X_1, \dots,X_k$ are in the kernel. In fact,
if $v\in E^{ss}$, by invariance of $d\omega_1,$ we have
$$
d\omega_1(X_i, v)=d\omega_i(D\phi^{ta_i}(X_i),D\phi^{ta_i}(v))=d\omega_1(X_i,D\phi^{ta_i}(v)).
$$
Hence,  the continuity of $d\omega_1$ and the fact that $a_i \in \mathcal{C}_0$,  imply that
$$
d\omega_1(X_i, v)=0.
$$
Similarly we obtain the same equality when $v\in E^{uu}$. Therefore, the claim is proved.

Finally, since the set of compact orbits of $\phi$  is dense in $M$ and $d\omega_1$ is continuous,
then it is sufficient to show that

{\it Claim 2: the 2-form $d\omega_1$ vanishes on each compact orbit.}
Let $p$ be a point in $M$ whose orbit is compact. Let $E=E_p^{ss}\oplus E_p^{uu}$,
$S=E_p^{ss}$, $U=E_p^{uu}$ and $a\in \mathcal{A}_0$ an element of the isotropy group of $p$.

Let $\ell$,   large enough, such that $f=(D\phi^v)^{\ell}$  satisfy all the conditions of Lemma 2.2
of \cite{ghys}.
If the 2-form $d\omega_1$ is non zero at $p$, then we must have that $|{\rm det} f|<1$, contradicting the fact
that $\phi^{ta}$ is volume preserving.
\end{proof}

 \section{Conclusion}
 \label{sec:conclusion}
 
In conclusion, we have the following sequence of implications for a given irreducible codimension one Anosov action of $\RR^k$:
 
 $$ 
\begin{array}{c}
     \mbox{$\phi$ is a linear suspension}\\
     \Downarrow\\
     \mbox{$\phi$ is transverse to a fibration over $\TT^k$} \\
     \Downarrow \\
 \mbox{$\phi$ is splitting} \\
 \Downarrow \\
 \mbox{$\phi$ is $1$-suspended}
\end{array}
 $$
 
Verjovsky conjecture states in particular that all these implications can be reversed. Let us discuss how it could be done for each of this.

\subsection{Are $1$-suspended actions splitting?}
This question reduces to the study of partially hyperbolic diffeomorphisms $f$ on a manifold $N$ whose central direction is tangent to a smooth locally free action of $\RR^l$ (with $l=k-1$) commuting with $f$. 

A first idea is that one could also expect to be able to prove that the $\RR^l$-action on $N$ is necessarily Anosov - but then one looses the codimension one property.

In the case $l=0$, ie. $k=1$, the positive answer comes essentially from the fact that  the stable foliation of an Anosov diffeomorphism has no holonomy and therefore, by a Theorem
of Sacksteder, admits an invariant transverse measure; hence is transversely affine, in a very peculiar way. 
The intermediate step for the proof of the existence of an invariant measure for (taut) codimension one foliations without holonomy is that the action of the fundamental group
on the leaf space is free. In the same line, it is easy to show that $1$-suspended actions  have the following remarkable property:

\begin{lemma}
Let $f_\ast: \cL^s \rightarrow \cL^s$ the diffeomorphism induced by any lifting of $f$ on the universal covering. Then, the fixed points of any positive iterate $f_\ast^p$ ($p >0$) are all 
repelling fixed points.
\end{lemma}

One could hope that this property, combined with other <<obvious>> properties of $f_\ast$, leads to the conclusion that $\cL^s$ is Hausdorff, ie. diffeomorphic to $\RR$. 
Then, every element of $\pi_1(N)$ would admit at most one fixed point. Indeed, if $x$ is a fixed point in $\cL^s$ of a non-trivial element $\gamma$ of $\pi_1(M)$, then it is also
a fixed point of $f_\ast^p$ with $p >0$ for some lifting of $f$, since non-planar orbits of $\RR^k$ are tori invariant by an iterate of $f$! Moreover, $f_\ast$ and $\gamma$ must commute, and thus, $x$ 
has to be the unique fixed point of $f_\ast$. We could then conclude, thanks to Solodov Theorem \ref{thm:solodov} that $\phi$ is splitting.

Observe that the case $k=2$, ie. $l=1$, is already interesting: our question then reduces essentially to partially hyperbolic diffeormorphisms $f$ of codimension one and whose central direction is smooth
and such that the leaves of the $1$-dimensional central foliation (which in this case exists necessarily) have non-trivial holonomy. 

\subsection{Are splitting Anosov actions of $\RR^k$ transverse to a fibration over $\TT^k$?}

When the action is splitting, we already know that the ambient manifold $M$ has a contractible universal covering (even homeomorphic to $\RR^n$), and a fundamental
group isomorphic to the semi-direct product of $\ZZ^n$ by $\ZZ^k$. It follows that $M$ is homotopy equivalent to the manifold suspension of the torus $\TT^n$ by a linear action
of $\ZZ^k$. By a foliation's classifying space argument, it can even be proved that the homotopy equivalence maps the stable foliation $\cF^s$ onto a codimension one linear
foliation. 

Moreover, as we have shown, we already know that the action is $1$-suspended. One line of proof could be to find sufficiently enough fibrations over the circle, transverse one to the other so that they altogether define a fibration over $\TT^k$ transverse to the action.

\subsection{Are actions transverse to a fibration topologically equivalent to a suspension of linear automorphisms?}

As we proved, a positive answer to this question would be given by finding an Anosov monodromy element (Lemma \ref{le:anosovmonodromy}).


\end{document}